\documentclass[american]{article}
\usepackage{babel,amssymb,amsfonts}
\usepackage{amsmath}
\usepackage{fancybox}
\usepackage{times}
\newtheorem{theorem}{Theorem}
\newtheorem{cor}{Corollary}
\newtheorem{prop}[cor]{Proposition}

\newcommand{\mr}{\mathrm}

\newcommand{\mb}{\mathbf}
\newcommand{\ms}{\mathsf}
\newcommand{\mf}{\mathfrak}

\newcommand{\bb}{\mathbb}
\newcommand{\ep}{\hfill$\square$}

\dateamerican
\title{Towards a nobabelian cohomology of forms}
\author{Mukul Patel\\ Department of Mathematics\\University of Georgia\\
Athens, GA  (U.S.A.)\\ patel@math.uga.edu\\ \date{\today}}
\begin{document}

\maketitle

\boldmath

\begin{abstract}

We consider a rather simple and natural coboundary operator, $\ms{d}$, on the Lie algebra valued differential forms on a manifold, which reduces to usual exterior derivative $d$ of such forms in the abelian case. Using the corresponding de Rham cohomology Lie superalgebra $\mr{H^*(M, \mr{G})}$, we obtain numerical \emph{smooth} invariants---as opposed to homotopy invariants---for manifolds. The corresponding Hodge theory yields finiteness of nonabelian Betti numbers. A genralized Poincar\'{e} lemma, along with a Poincar\'{e} duality, a Mayer-Vietoris, and a particularly empowered Bockstein makes our cohomology \emph{computable}.  Bockstein also allows us to relate (nonabelian) diffeomorphism invariants to (abelian) homotopy invariants.

\end{abstract}

\setcounter{secnumdepth}{1}
\setcounter{tocdepth}{3}

%\tableofcontents

\section{Introduction}
Throughout this paper, $\mr{G}$ will denote a real Lie group with Lie algebra $\mf{g}.$
Then, the graded space of  $\mathfrak{g}$-valued differential forms on a smooth manifold $\mr{M}$ is a $\bb{Z}$-graded Lie superalgebra with respect to a natural bracket operation. (See \ref{2}, below for definition.) Indeed, this algebra is the tensor product of the exterior algebra of usual $\mathbb{R}$-valued forms with $\mathfrak{g}$, and the usual exterior derivative operator is just $d=$ \math${d}
$ \otimes$ $\mr{Id_{\mr{G}}}$, where \math${d} is the usual exterior derivative of $\bb{R}$-valued forms .  Since the differential $d$ does not involve
the Lie algebra structure of $\mathfrak{g}$, the cohomology Lie superalgebra based on it is just the tensor product of the usual de Rham cohomology algebra with $\mathfrak{g}$. Obviously, this can not distinguish between two
manifolds any more than the de Rham cohomology.

In this paper, we consider another natural differential $\ms{d}$, which directly involves the Lie algebra structure of $\mathfrak{g}$, and reduces to $d=$ \math${d}
$ \otimes$ $\mr{Id_{\mf{g}}}$ when $\mathfrak{g}$ is abelian, and to  the usual \math${d} when $\mf{g}= \bb{R}.$ Thus, we obtain an essentially new cohomology Lie superalgebra which is a smooth invariant of the manifold $\mr{M}$. We note, however, that this is not a homotopy invariant in general when $\mr{G},$ and hence $\mf{g},$ is nonabelian. Indeed, \textbf{this cohomology is homotopy invariant if and only if $\mr{G}$ is abelian.} See the last section for a discussion of this phenomenon.

Once we have identified this nonabelian coboundary operator, most of the results in this paper are routine verifications following the
proofs of the corresponding statements about the usual de Rham cohomology.

\section{Nonabelian cohomology of forms} \label{2}

A $\bb{Z}$-\textbf{graded Lie superalgebra} over a field $\mr{K}$, is a $\bb{Z}$-graded vector $\mr{K}$-vector space $\mr{L}= \mr{\sum _{p=0}^{n} }\mr{L}_{\mr{p}}$, along with a (non-associative) product 
$[\hspace{2pt}, \hspace{2pt}]$, such that for $a\in\mr{L}_{p}$ and $ b\in\mr{L}_{q},$ we have $[a,b]\in\mr{L}_{p+q}$, and:
\begin{enumerate}
\item $[a, b] =  -(-1)^{pq}\hspace{1pt}[b, a]$\hfill (super-commutation)
\item $[a,[b, c]] = [[a, b], c] +  (-1)^{pq}\hspace{1pt}[b, [a, c]].$\hfill (Jacobi Identity)
\end{enumerate}
(See \cite{kac,cns} for details on Lie superalgebras.)

The graded space $\mr{\Lambda^{*}(M, \mf{g})},$ of $\mf{g}$-valued forms on $\mr{M}$ is a 
$\bb{Z}$-graded Lie superalgebra, where the bracket is defined as the composition:
\[[\hspace{2pt}, \hspace{2pt}]:  \Lambda^{\mr{p}} \mr{(M, \mf{g}) \times \Lambda^{q} (M, \mf{g})}\overset{\Lambda}{\longrightarrow}\mr{\Lambda ^{p+q}(M, \mf{g} \otimes \mf{g}})\overset{[\hspace{2pt},\hspace{2pt}]}{\longrightarrow }\Lambda^{p+q}(\mr{M}, \mf{g})\]
(See \cite{bleecker} for details on this algebra.)

As mentioned in the introduction, this algebra is simply the tensor product of the real-valued exterior algebra with the Lie algebra $\mf{g}$, which has the trivially extended differential $d=$ \math${d}
$ \otimes$ $\mr{Id_{\mr{G}}}$.
To define our nonabelian differential, we start with $\mr{G}$-equivariant $\mf{g}$-valued forms on the trivial principal bundle $\mr{M\times G}$, such that if any of the vector fields these forms act on is vertical, then the value of the form is vanishes.  We denote the algebra of these forms by $\mr{\bar{\Lambda}^{*}(M\times G, \mf{g})}$. Note that pulling back forms along the projection $\pi : \mr{M \times G \rightarrow M}$ induces a natural isomorphism:
\[ \mr{\Lambda^{*}(M, \mf{g})} \cong  \mr{\bar{\Lambda}^{*}(M\times G, \mf{g})}:
\alpha \mapsto \bar{\alpha}= \pi^*{\alpha}.
\]
Now onwards, we will identify these algebras, and think of each $\alpha \in  \mr{\Lambda^{*}(M, \mf{g})}$  as being the form $\pi^*{\alpha}$ living on the principal bundle. To emphasize this point we will write  $\mr{\Lambda^{*}(M, \mr{G})}$ for any of these two algerbas.

Let $\ms{d}\alpha \in \mr{\Lambda^{k}(M, \mr{G})}$ be the covariant exterior derivative
of $\alpha$ with respect to the canonical flat connection on the principal bundle  $\mr{M\times G}$,  given by the $\mf{g}$-valued connection form $\theta$.  (See \cite{kn}.)
This is the differential we will use to define a cohomology.

\begin{prop}
 \[\ms{d}\alpha =  d\alpha +  [\theta, \alpha],\]
and, consequently,
\[\ms{d}^{2}\alpha = [\Theta, \alpha]= 0,\]
\end{prop}
\textbf{Proof:}
The first assertion obtains by evaluating both sides on different combinations of verticle and horizontal forms. The first equality in the second assertion follows by direct calculation
and recalling that $\Theta = d\mr{\theta + \frac{1}{2}[\theta, \theta].}$ The second equality of the second assertion holds because the canonical connection $\theta$, is flat, \textit{i.e.,} its curvature $\Theta $ vanishes.\hfill$\square$

It is immediate from the above formula, that 
\begin{prop}
The operator $\ms{d}$ is an (anti) derivation of order $1$, and,
\[\ms{d}^{2} = 0.\]\hfill$\square$
\end{prop}

Note that for abelian $\mf{g}$, we have $\ms{d} = d.$  In particular, when $\mr{G} = \mf{g}= \bb{R},$ then $\ms{d}$ is the usual
de Rham coboundary operator on $\mr{M}.$

\begin{prop}
Let $\mr{f: M \rightarrow N}$ be a smooth map of smooth manifolds. Then,
\[f^{*}\ms{d} = \ms{d}f^{*},\]
where $\mr\mb{{f^*:\Lambda^{*}(N, \mr{G}) \rightarrow \Lambda^{*}(M, \mr{G})}}$
is given by the pullback of forms along $\mr{f}.$
\end{prop}
\textbf{Proof:}
This is immediate because the canonical connection $\theta$ on $\mr{M\times G}$ is just the pullback of the Maurer-Cartan form on $\mr{G}$ along the projection $\mr{M\times G\rightarrow G}$. \hfill$\square$

The \textbf{cohomology algebra of $\mr{\mb{M}}$ with values in $\mr{G}$} is the cohomolology Lie superalgebra, $\mr{\mb{H^*(M, \mr{G})}}$, of the cochain complex $\mr{\mb{(\Lambda^{*}(M, \mr{G})}, \ms{d})}$. We will call this cohomology theory the \textbf{$\mr{G}$-cohomology theory.}

\begin{prop}
For a fixed $\mr{G}$, the functor
$\mr{M \mapsto H^*(M, \mr{G})}$ is a contravariant functor from the category of smooth manifolds to the category of Lie superalgebras. For a fixed $\mr{M,}$ the functor $ \mr{G \mapsto  H^*(M, \mr{G})}$ is a covariant functor from the category of Lie groups to the category of Lie superalgebras.  

\end{prop}
\textbf{Proof:}
Define the cohomology maps by pullbacks of forms, and apply the preceding proposition.\ep

The Poincar\'{e} lemma below shows that for nonabelian $\mr{G}$, the cohomology is not homotopy invariant, but only diffeomorphism invariant,  while the cohomology with abelian $\mr{G}$ is homotopy invariant. Thus, maps between Lie groups, one of which is abelian, gives relationships between homotopy invariants and diffeomorphism invariants.

Now we turn to the usual staple of facts that aid computation of the cohomology algerbas.

\section{Aids to computation}
\subsubsection{Hodge theory}\label{hodge}
Note that $\mr{(\Lambda^{*}(M, \mr{G}), \ms{d})}$ is an elliptic complex, so that classical Hodge theory can be immediately extended to it \cite{warner}, \cite{wells}. We briefly recount the main points in our context.

Let us assume that $\mr{M}$ is oriented, with a  Riemannian metric, and that $\mr{G}$ is semisimple compact. Then,  the Cartan-Killing form on $\mf{g}$ is nondegenerate and (negative) definite, which allows us to define an inner product on $ \Lambda^{*}\mr{(M, \mr{G})}$ such that its homogeneous components are orthogonal. We can also define the Hodge $*$-operator: 
\[*: \mr{\Lambda^{k}(M, \mr{G}) \rightarrow \Lambda^{n-k}(M, \mr{G})}\] for each $\mr{0\leq k\leq n}$ in the usual fashion; see \cite{warner}. Now define \[\mr{\ms{d}^{*} : \Lambda^{k}(M, \mr{G}) \rightarrow \Lambda^{k-1}(M, \mr{G})}\] by $ \ms{d}^{*}\alpha = (-1)^{n(p+1)+1}(* \ms{d} *)(\alpha)$. This operator is the formal adjoint of $\ms{d}$. Finally, for the \textbf{Laplace-Beltrami operator,}
$\triangle := \ms{d}^{*}\ms{d}+\ms{d}\hspace{1pt}\ms{d}^{*},$ the generalized Hodge theorem asserts that \[\mr{Ker}\triangle \cong \mr{H^*(M, \mr{G})}: \alpha\mapsto[\alpha]\] where $[\alpha]$ is the cohomology class of $\alpha$. See \cite{wells}.
Note that the left hand side of the isomorphism is the kernel of an elliptic operator, and as such, is finite dimensional. Consequently, each cohomology space $\mr{H^k(M, \mr{G})}$ is finite dimensional.  Using the Bockstein and Mayer-Vietoris long exact sequences, this finiteness result can be extended to the general case where $\mr{G}$ is not assumed to be semisimple compact . The forms in the kernel of Laplace-Beltrami operator are called \textbf{harmonic forms}. We will need to return to these matters in defining the Poincar'{e} duality.

\subsubsection{Poincar\'{e} lemma}
Let $\mr{c}\mf{g}$ be the center of a Lie algebra $\mf{g}$, and $\mf{[g, g]'}$ be the subalgebra generated by all the elements of $\mf{g}$ that commute with the commutator ideal $\mf{[g, g]}$.  Then, a straightforward computation shows that the following version of Poincar\'{e} lemma holds.
\begin{theorem}[Poincar\'{e} lemma]\label{pl}
For any nonnegative integer $\mr{n},$ 
\[\mr{H^k(\bb{R}^{\mr{n}}, \mr{G})} = \left \{ \begin{array}{ll}
\mr{c}\mf{g }&\hspace{5pt} \mbox{if \hspace{5pt}$\mr{k = 0}$} \\
\mf{[g, g]'}/\mr{c}\mf{g}&\hspace{5pt} \mbox{if \hspace{5pt}$\mr{0<k<n}$}\\
0  &\hspace{5pt}\textit{otherwise.}
\end{array}
\right .\]

\end{theorem}

\ep

Thus, it follows that $\mr{H^k(M, \mr{G})}$ are not homotopy invariants unless $\mr{G}$ is an abelian Lie group.
Now we define  $k$-th \textbf{$\mr{G}$-Betti number} of $\mr{M}$,
\[\mr{b_{k}^{\mr{G}} :=  dim\hspace{2pt}H^k(M, \mr{G}),}\]
and the $\mr{G}$-\textbf{Euler Characteristic} of $\mr{M}$,
\[\mr{\chi(M, \mr{G}) := \sum(-1)^{k} b_{k}^{\mr{G}}}.\]

Note that $\mr{b_{k}^{\bb{R}}= b_{k}^{}},$ the usual Betti numbers, and
$ \mr{\chi(M, \bb{R})= \chi(M),}$ the Euler chracteristic of $\mr{M}.$
We record the following information contained in Poincar\'{e} lemma:

\begin{theorem}\label{betti}
Set dimension of $\mf{g} = \mr{g}$ and dimension of $\mr{c}\mf{g} = \mr{g}_{\mr{c}}$, and dimension of $\mf{[g, g]'}= \mr{g'}.$
Then, the $\mr{G}$-Betti numbers and $\mr{G}$-Euler characteristic of $\bb{R}^{\mr{n}}$ are given by
\[\mr{b^{\mr{G}}_{k}(\bb{R}^{\mr{n}})} = \left \{ \begin{array}{ll}
\mr{g_{c}}&\hspace{5pt} \mbox{if \hspace{5pt}$\mr{k = 0}$} \\
\mr{g'-g_{c}} &\hspace{5pt} \mbox{if \hspace{5pt}$\mr{0<k<n}$} \\
0  &\hspace{5pt}\textit{otherwise.}
\end{array}
\right .\]
and
\[\mr{\chi(\bb{R}^{\mr{n}}, \mr{G})} = \left \{ \begin{array}{ll}
\mr{g_{c}}&\hspace{5pt} \mbox{if \hspace{5pt}$\mr{ n = 0}$} \\
\mr{2g_{c}-g'} &\hspace{5pt} \mbox{if \hspace{5pt}$\mr{n>0 }$ \textit{is even}}\\
\mr{g}_{c}&\hspace{5pt} \mbox{if \hspace{5pt}$\mr{n>0 }$ \textit{is odd}}
\end{array}
\right .\]\ep

\end{theorem}

Thus, we obtain numerical diffeomorphism invariants.  For an \emph{abelian} $\mr{G}$, with dimension $\mr{g},$ we have
$\mr{b_{k}^{\mr{G}} =  g.b_{\mr{k}}^{}}$, making the cohomology spaces and the $\mr{G}$-Euler Characteristic homotopy invariants.

\subsubsection{Cohomology exact sequences}

Both the Mayer-Vietoris and Bockstein long exact sequences continue to be exact in our general setting.
However, a word on Bockstein is in order\hspace{1pt}:
In the defnition of abelian cohomology functor, the second variable is only nominally involved (in case of torsionless coefficients). Since the coefficients are nontrivially involved in the nonabelian case, the Bockstein can be expected to yield substantially more information on the cohomology functor. In particular, the second variable (the coefficient algebra $\mf{g}$) is actively involved in cohomology operations, i.e. natural transformations on the functor. As a quick application of Bockstein, we note that in view of the comments following Theorem \ref{betti}, the Bockstein for the canonical exact sequence 
\[0\longrightarrow\mr{[G,G]}\longrightarrow \mr{G}\longrightarrow\mr{G}/[\mr{G},\mr{G}]\longrightarrow 0\] 
or, more generally, for any exact sequence of Lie groups, one of which is abelian, aids in computing relationships between diffeomorphism invariants and homotopy invariants, \textit{e.g.} the corresponding Betti numbers and Euler characteristic.

\subsubsection{Poincar\'{e} duality}
For an oreinted compact $\mr{n}$-manifold $\mr{M}$ without boundary, Poincar\'{e} duality is given by a nonsingular pairing 
\[\mr{H^k(M, \mr{G})\otimes H^{n-k}(M, \mr{G})}\longrightarrow \bb{R},\]
which gives an isomorphism
\[\mr{H^k(M, \mr{G})\cong Hom(H^{n-k}(M, \mr{G}), \bb{R})}.\]
The pairing is defined as in the de Rham cohomology:
Let $\alpha$ be the unique harmonic form representing a cohomology class $[\alpha]$ in $\mr{H^k(M, \mr{G})}$,
and let $\beta$ be the unique harmonic form representing a cohomology class $[\beta]$ in
$\mr{H^{n-k}(M, \mr{G})}$. Then we define the duality pairing by,
\[([\alpha], [\beta])\, \mapsto \;<\!\alpha , \beta^*\!\!> \]
Where the $< , >$ is the inner product defined above in the section on Hodge theory.
For $\alpha \neq 0$,
\[([\alpha], [\alpha^*])\: \mapsto \; <\!\!\alpha, \alpha\!\!> \;= \|\alpha\| ^{2}\neq  0.\]
Thus, the pairing is nonsingular. There is, however, a more subtle duality given by a pairing
\[\mr{H^k(M, \mr{G})\otimes H^{n-k}(M, \mr{G})}\longrightarrow \mf{g},\]
where the cohomology spaces are considered as $\mf{g}$-modules. The use of this duality involves detailed knowledge of $\mf{g}$-modules for the specific Lie algebra $\mf{g}$. 

%\subsubsection{Leray-Hirsch spectral sequence and Kunneth formula}

%\subsubsection{Atiyah-Hirzebruch spectral sequence}

\section{New invariants for manifolds}

We have already looked at some numerical invariants related to individual homogeneous components of the cohomology  Lie superalgebra. To those, we add the invariants of the whole cohomology algebra. We note that it is a $\bb{Z}$-graded Lie superalgebra,  the even and odd parts of which are the direct sum of its even and odd $\bb{Z}$-graded components respectively:
\[\mr{H^*(M, \mr{G}) =  H^{even}(M, \mr{G}) \oplus H^{odd}(M, \mr{G}).}\]
The even component is a Lie algebra under the bracket product of the Lie superalgebra. The superalgerba structure gives a canonical representation of this ($2\bb{Z}$-graded)  Lie algebra on the odd part, and the bracket of any two elements of the odd part is in an even part. Thus, we have several new invariants of manifolds, including the invariants of the Lie algerba constituting the even part. Notice that the cohomology algebra is clearly \emph{nilpotent,} and so is its even component (as a Lie algebra). Recent progress in the classification of nilpotent Lie algebras (\cite{gk}) can be fruitfully used in this context.

Since the cohomology theories with abelian  $\mf{g}$ are related by a universal coefficient theorem, they contain roughly the same information. However, in our general nonabelian setting, the Poincar\'{e} lemma shows that there is no universal coefficient theorem. This is a good indication that, in general, different choices of the nonabelian coefficents $\mr{G}$ yield essentially different cohomology theories, although related by functoriality in the second variable of the bifunctor $\mr{H^*(M, \mr{G}})$. 
Thus, it seems desirable to determine the extent to which this functor encodes purely the information on smooth structure of $\mr{M},$ and to somehow `mod out' the purely algerbaic information coming from the group $\mr{G}.$ One way to do this is to consider $\mr{H^*(M, \mr{G})}$ as a functor from
the category $\ms{M}$ of smooth manifolds to a functor category,
\[\mr{H : \ms{M}\rightarrow \ms{S}:\;\;M \mapsto HM},\]
the objects of $\ms{S}$ being functors  from the category $\ms{G}$ of Lie groups to the category $\ms{SG}$ of Lie superlagebras, where $\mr{HM}$ is the functor
\[\mr{HM : \ms{G}\rightarrow \ms{SG}: \;\;\mr{G} \mapsto  H^*(M, \mr{G})}.\]
Furthermore, the functor $\mr{HM}$ yields, via Bockstein, a homomorphism  from the\linebreak
Grothendieck group of $\ms{G}$ to the Grothendieck group of the category $\ms{SG} :$
\[ \mr{HM : K(\ms{G})\rightarrow K(\ms{SG})}.\]
This homomorphism is a smooth invariant of $\mr{M}$.  Thus, we introduce the functor
\[\mr{M \mapsto HM},\]
the right hand side being the homomorphism just introduced. Now, Mayer-Vietoris and Bockstein for $\mr{H^*(M, \mr{G})}$ yield a Mayer-Vietoris for this functor.

\section{General remarks on nonabelian cohomology theories}

An important point is that the nonabelian character of the coeffcients is manifested in the muliplicative aspect of the \emph{whole} cohomology algebra, which, in the abelian case, decomposes into it's graded component spaces. This indicates
that the real object of a (possibly nonabelian) cohomology theory is the entire cohomology algerba and that it is fortuitous that in the abelian case each cohomology set has an algebraic structure. 

This sits well with the fact that the first cohomology set with coeffcients in a nonabelian (topological) group (or a sheaf of nonabelian groups) has no algebraic structure at all, which forces us to look for an algebraic structure on the union of all prospective higher cohomology sets.  

Recall that the cohomology set $\mr{H^1(X, G)}$ of a space $\mr{X}$ classifies principal $\mr{G}$-bundles over $\mr{X}$, and is represented by classifying space $\mr{BG}:  \mr{H^1(X, G) \cong [X, BG]}$. Here, $\mr{G}$ is a locally compact topological group. When this is an abelian group, $\mr{BG}$ is (an abelian) group, and its classifying space $\mr{B^2G := B(BG)}$ represents the second cohomology group $\mr{H^2(X, G) \cong [X, B^2G]}$. Now, $\mr{B^2G}$ is an abelian group, and this process can be continued to yield representations of the higher cohomomology functors.  However, for nonabelian $\mr{G}$, the space $\mr{BG}$ is not even an $\mr{H}$-space, let alone a group. This circumstance abruptly truncates the construction of higher ``classifying spaces"  and is reflected in the fact that there is no satisfactory algebraic theory of higher nonabelian cohomology.  In \cite{me3}, a construction of natural higher classifying spaces is presented, and then is used to \emph{define} the corresponding higher cohomology sets. The latter classify higher torsors---analogous $\mr{H^1(X, G)}$ classifying principal bundles. Needless to say that none of these cohomology sets have any algebraic structure (except that they are pointed sets). This indicates that if there is any hope of finding a usable nonabelian cohomology theory, we must look for a classifying object which would represent (union of) \emph{all} the cohomology sets at once. In order for this structure to be any more useful than abelian cohomologies, we need noncommutative multiplicative strucure. 

A tentative start towards a general nonabelian theory along these lines is presented in forthcoming papers \cite{me,me2}.


\begin {thebibliography}{BB}

\bibitem[B]{bleecker} Bleecker, D., 1981, \textit{Gauge theory and variational principiles,} Addison-Wesley, Boston.
\bibitem[CNS]{cns} Corwin, L., Ne'eman, Y. and S. Sternberg, 1975, Graded Lie algebras in mathematics and physics, \textit{Reviews of Modern Physics,} 47, No. 3, 573-603
\bibitem[GK]{gk} Goze, M. and Y. Khakimdjanov, 1996, \textit{Nilpotent Algebras,} Kluwer Academic, Boston.
\bibitem[K]{kac} Kac, V., 1977, Lie Superalgebras \textit{Advnaces in Math.} 26, No 1, 8-96
\bibitem[KN]{kn} Kobayashi, S., and K. Nomizu, 1991, \textit{Foundations of differential geometry - I,} John Wiley, New York.
\bibitem[M1]{me}  Patel, M., \textit{Nonabelian \v{C}ech cohomology,} (To appear)
\bibitem[M2]{me2} Patel, M., \textit{On nonabelian cohomology,} (To appear)
\bibitem[M3]{me3} Patel, M., \textit{Higher torsors, nonabelian cohomology, and classifying spaces,}(To appear)
\bibitem[W]{warner} Warner, F. W., 1983, \textit{Foundations of differentiable manifolds and Lie groups,} Springer Verlag, New York.
\bibitem[WL]{wells} Wells, R. O., 1980, \textit{Differential analysis on complex manifolds,} Springer Verlag, New York.


\end{thebibliography}
\end{document}